\documentclass[12pt,leqno]{article} % {mf}
\usepackage{amssymb}
\usepackage{amsmath}
\usepackage{mathrsfs}
\newcommand{\R}{{\bf R}}

\newcommand{\mr}{\mathbb{R}}
\newcommand{\cf}{{\cal F}}
\newcommand{\ce}{{\cal E}}
\newcommand{\cR}{{\cal R}}
\newcommand{\ca}{{\cal A}}

\newcommand{\ma}{\mathscr{A}}
\newcommand{\md}{\mathscr{D}}

\numberwithin{equation}{section}
\newtheorem{theo}{\sc Theorem}[section]

\newtheorem{coro}[theo]{\sc Corollary}%[section]
\newtheorem{assump}[theo]{\sc Assumption}%[section]
\newtheorem{remark}[theo]{\sc Remark}%[section]
\newtheorem{ex}[theo]{\sc Example}%[section]
\newtheorem{defi}[theo]{\sc Definition}%[section]
\newtheorem{prob}[theo]{\sc Problem}%[section]
\newcommand{\eof}{\hfill $\Box$ \par}

\allowdisplaybreaks

\def\mE{\mathscr{E}}

\begin{document}
\title{A stochastic maximum principle via Malliavin calculus}
%{A Malliavin calculus approach to a general maximum principle for stochastic control}
\author{Thilo Meyer-Brandis$^{1)}$
%\thanks{Center of Mathematics for Applications (CMA), University of
%Oslo, Box 1053 Blindern, N-0316 Oslo, Norway. Email: $<$meyerbr@math.uio.no$>$.}
\and
Bernt {\O}ksendal$^{1),2)}$%\thanks{Center of Mathematics for Applications (CMA), University of
%Oslo, Box 1053 Blindern, N-0316 Oslo, Norway. Email: $<$oksendal@math.uio.no$>$.}
%\thanks{Norwegian School of Economics and Business Administration (NHH), Helleveien 30, N-5045
%Bergen, Norway}
\and
Xun Yu Zhou$^{3)}$}
%\thanks{Mathematical Institute, University of Oxford,
%24-29 St Giles', Oxford OX1 3LB, UK, and
%Dept of Systems Engineering and Engineering Management, The Chinese
%University of Hong Kong, Shatin, Hong Kong.  Email: $<$zhouxy@maths.ox.ac.uk$>$. }}
\date{November 19, 2009}%{4 July 2007}

\footnotetext{$^{1)}$Center of Mathematics for Applications (CMA), University of
Oslo, Box 1053 Blindern, N-0316 Oslo, Norway. Email: $<$meyerbr@math.uio.no$>$,
$<$oksendal@math.uio.no$>$.}

\footnotetext{$^{2)}$Norwegian School of Economics and Business Administration (NHH), Helleveien 30, N-5045 Bergen, Norway.}

\footnotetext{$^{3)}$Mathematical Institute, University of Oxford,
24-29 St Giles', Oxford OX1 3LB, UK, and
Dept of Systems Engineering and Engineering Management, The Chinese
University of Hong Kong, Shatin, Hong Kong.  Email: $<$zhouxy@maths.ox.ac.uk$>$. }

\maketitle
\begin{abstract}
This paper considers a controlled It\^o-L\'evy process where the
information available to the controller is possibly less than the
overall information. All the system coefficients and the objective performance
functional are allowed to be random, possibly non-Markovian.
Malliavin calculus is employed to derive a  maximum
principle for the optimal control of such a system where the adjoint process is explicitly expressed.
\medskip

\noindent {\bf Mathematics Subject Classification 2000}: 93E20, 60H10, 60HXX, 60J75
\medskip

\noindent{\bf Key words}: Malliavin calculus, maximum principle, stochastic control, jump diffusion, partial information
\end{abstract}
%%%%%%%%\iffalse
\section{Introduction}
Suppose the state process $X(t)=X^{(u)}(t,\omega)$;  $t\geq 0$, $\omega\in \Omega$, is a controlled It\^{o}-L\'{e}vy process in $\mathbb{R}$ of the form
\begin{equation}\label{fomu1}
\left\{\begin{array}{l}
dX(t)=b(t,X(t),u(t),\omega)dt+\sigma(t,X(t),u(t),\omega)dB(t)\\
\hspace{2cm}+\int_{\mr_0} \theta(t,X(t^-),u(t^-),z,\omega)\tilde{N}(dt,dz);\\
X(0)=x\in \mr.\end{array}
\right.
\end{equation}
Here $\mr_0=\mr-\{0\}$, $B(t)=B(t,\omega)$, and $\eta(t)=\eta(t,\omega)$, given by
\begin{eqnarray}
\eta(t)=\int^t_0\int_{\mr_0}z\tilde{N}(ds,dz);\; t\geq 0,\;\omega\in\Omega,
\end{eqnarray}
are a 1-dimensional Brownian motion and an independent pure jump L\'{e}vy martingale,
respectively, on a given filtered probability space $(\Omega, {\cal F}, \{
{\cf}_t \}_{t\geq 0},P).$ Thus
\begin{equation}
\tilde{N}(dt,dz):=N(dt,dz)-\nu(dz)dt
\end{equation}
is the {\em compensated jump measure\/} of $\eta(\cdot)$, where $N(dt,dz)$ is the
{\em jump measure\/} and $\nu(dz)$ the {\em L\'{e}vy measure\/} of the L\'{e}vy
process $\eta(\cdot).$ The process $u(t)$ is our control process, assumed to be ${\cal F}_t$-adapted  and have values in a given open convex set $U\subset \mr.$ The coefficients $b:[0,T]\times \mr\times U\times \Omega\rightarrow \mr$, $\sigma:[0,T]\times\mr\times U\times \Omega\rightarrow \mr$ and $\theta:[0,T]\times\mr\times U\times \mr_0\times \Omega$ are given ${\cal F}_t$-predictable processes.

We refer to [YZ] and [\O{}S] for more information about stochastic control of It\^{o} diffusions
and jump diffusions, respectively.
Let $T>0$ be a given constant. For simplicity, we assume that
\begin{eqnarray}
\int_{\mr_0} z^2\nu(dz)<\infty.
\end{eqnarray}
Suppose in addition that we are given a subfiltration
$$
\mE_t\subseteq {\cf}_t,\qquad t\in[0,T]
$$
representing the information available to the controller at time $t$ and satisfying the usual conditions. For example, we could have
$$
\mE_t={\cf}_{(t-\delta)^+};\quad  t\in[0,T],\;\delta>0 \mbox{ is a constant},
$$
meaning that the controller gets a delayed information compared to ${\cf}_t$.

Let
 $\mathscr{A}=\ma_\mE$
 denote a given family of controls, contained in the
set of $\mE_t$-adapted c\`adl\`ag controls $u(\cdot)$ such that (\ref{fomu1}) has a unique strong solution up to time $T$.
Suppose we are given a performance functional of the form
\begin{eqnarray}
J(u)=E\Big[\int^T_0f(t,X(t),u(t),\omega)dt+g(X(T),\omega)\Big];\qquad u\in \ma_\mE,
\end{eqnarray}
where $E=E_P$ denotes expectation with respect to $P$ and
$f:[0,T]\times \mr\times U\times \Omega\rightarrow \mr$ and $g:\mr\times\Omega\rightarrow\mr$
are given ${\cf}_t-$adapted processes with
$$
E\Big[\int^T_0|f(t,X(t),u(t))|dt+|g(X(T))|\Big]<\infty\qquad \mbox{for all}\quad u\in\ma_\mE.
$$

The partial information control problem we consider is the following:
\begin{prob}%PROBLEM
\label{prob1}
Find $\Phi_\mE \in \mr$ and $u^*\in \ma_\mE$ (if it exists)
such that
\begin{equation}
\Phi_\mE=\sup_{u\in\ma_\mE}J(u)=J(u^*).
\end{equation}
\end{prob}

Note that since we allow $b,\sigma,\theta,f$ and $g$ to be stochastic processes and also because our controls must be $\mE_t-$adapted, this problem is not of Markovian type and hence cannot be solved by dynamic programming.
We instead investigate the maximum principle, and derive an explicit form for the adjoint process. The approach we employ is Malliavin calculus which enables us to express the duality involved via the Malliavin derivative.
Our paper is related to the recent paper [B\O], where a maximum principle for partial information control is obtained. However, that paper assumes the existence of a solution of the adjoint equations. This is an assumption which often fails in the partial information case.

We emphasize that our problem should be distinguished from the {\em partial observation\/} control problems, where it is assumed that the controls are based on noisy observation of the state process. For the latter type of problems, there is a rich literature see e.g.  [BEK], [B], [KX], [L], [PQ], [T]. Note that the methods and results in the partial observation case do not apply to our situation. On the other hand, there are several existing works on stochastic maximum principle (either completely or partially observed) where adjoint processes are explicitly expressed [BEK], [EK], [L], [T]. However, these works all essentially employ stochastic flows technique, over which the Malliavin calculus has an advantage in terms of numerical computations (see, e.g., [FLLL]).

% Our paper is related to the recent paper [B\O], where a maximum principle for partial information control is obtained. However, that paper assumes the existence of a solution of the adjoint equations. This is an assumption which often fails in the partial information case.
%
% In this paper we use Malliavin calculus to obtain  a  maximum principle for this general non-Markovian partial information stochastic control problem. See Theorem \ref{GPIMPI}.

\section{A brief review of Malliavin calculus for L\'{e}vy processes}

In this section we recall the basic definition and properties of Malliavin calculus for L\'{e}vy
processes related to this paper, for reader's convenience.

In view of the L\'{e}vy--It\^{o} decomposition theorem, which states that any L\'{e}vy process
$Y(t)$ with
$$
E[Y^2(t)]<\infty \quad \mbox{for all}\quad t
$$
can be written
$$
Y(t)=at+bB(t)+\int^t_0\int_{\mr_0}z\tilde{N}(ds,dz)
$$
with constants $a$ and $b$,
we see that it suffices to deal with Malliavin calculus for $B(\cdot)$ and for
$$
\eta(\cdot):=\int_0\int_{\mr_0}z\tilde{N}(ds,dz)
$$
separately.

A general reference for this presentation is [N], [BDL{\O}P] and [DM{\O}P]. See also the forthcoming
book [D{\O}P].

\subsection{Malliavin calculus for $B(\cdot)$}

A natural starting point is the Wiener-It\^{o} chaos expansion theorem, which states that any
$F\in L^2({\cal F}_T,P)$ can be written
\begin{eqnarray}
F=\sum_{n=0}^{\infty}I_n(f_n)
\end{eqnarray}
for a unique sequence of symmetric deterministic functions $f_n\in L^2(\lambda^n)$,
where $\lambda$ is Lebesgue measure on $[0,T]$ and
\begin{eqnarray}
I_n(f_n)=n!\int^T_0\int^{t_n}_0\cdots\int^{t_2}_0f_n(t_1,\cdots,t_n)dB(t_1)dB(t_2)\cdots dB(t_n)
\end{eqnarray}
(the $n$-times iterated integral of $f_n$ with respect to $B(\cdot)$) for $n=1,2,\ldots$
and $I_0(f_0)=f_0$ when $f_0$ is a constant.

Moreover, we have the isometry
\begin{equation}
E[F^2]=||F||^2_{L^2(p)}=\sum^\infty_{n=0}n!||f_n||^2_{L^2(\lambda^n)}.
\end{equation}

\begin{defi}[Malliavin derivative $D_t$]  %{Definition 2.1}
\hfill\break
{\rm Let
%$\mathscr{D}_{1,2}=\md^{(B)}_{1,2}$
$\md^{(B)}_{1,2}$ be the space of all $F\in L^2({\cf}_T,P)$
such that its chaos expansion (2.1) satisfies
\begin{eqnarray}
||F||^2_{\md^{(B)}_{1,2}}:=\sum^\infty_{n=1}n n!||f_n||^2_{L^2(\lambda^n)}<\infty.
\end{eqnarray}

For $F\in \md^{(B)}_{1,2}$ and $t\in [0,T]$, we define the {\em Malliavin derivative\/} of F at t (with respect to $B(\cdot)$),
$D_tF,$ by
\begin{eqnarray}\label{Malder}
D_tF=\sum^\infty_{n=1}nI_{n-1}(f_n(\cdot,t)),
\end{eqnarray}
where the notation $I_{n-1}(f_n(\cdot,t))$ means that we apply the $(n-1)$-times iterated
integral to the first $n-1$ variables $t_1,\cdots, t_{n-1}$ of $f_n(t_1,t_2,\cdots,t_n)$
and keep the last variable $t_n=t$ as a parameter.}
\end{defi}

One can easily check that
\begin{eqnarray}\label{isometry}
E\Big[\int^T_0(D_tF)^2dt\Big]=\sum^\infty_{n=1}n n!||f_n||^2_{L^2(\lambda^n)}=||F||^2_{\md^{(B)}_{1,2}},
\end{eqnarray}
so $(t,\omega)\rightarrow D_tF(\omega)$ belongs to $L^2(\lambda \times P)$.

%\underline{Example 2.2}
\begin{ex}
\rm
If $F=\int^T_0f(t)dB(t)$ with $f\in L^2(\lambda)$ deterministic, then
$$D_t F=f(t) \mbox{ for } a.a. \,t\in[0,T].$$
More generally, if $u(s)$ is Skorohod integrable, $u(s)\in \md_{1,2}$ for $a.a. \; s$ and $D_tu(s)$
is Skorohod integrable for $a.a. \;t$, then
\begin{equation}
D_t\Big(\int_0^Tu(s)\delta B(s)\Big)=\int_0^TD_tu(s)\delta B(s)+u(t)\;
\mbox{for a.a. $(t,\omega)$},
\end{equation}
where $\int_0^T\psi(s)\delta B(s)$ denotes the Skorohod integral of $\psi$
with respect to $B(\cdot)$. (See [N], page 35--38 for a definition of Skorohod integrals and for more
details.)
\end{ex}

Some other basic properties of the Malliavin derivative $D_t$ are the following:
\begin{enumerate}
\item [(i)] {\bf Chain rule } ([N], page 29)\\
Suppose $F_1,\ldots,F_m\in\md^{(B)}_{1,2}$ and that $\psi:\mr^m\rightarrow \mr$ is $C^1$ %$e^\bot$
with bounded partial derivatives. Then $\psi(F_1,\cdots, F_m) \in \md_{1,2}$ and
\begin{eqnarray}
D_t\psi(F_1,\cdots, F_m)=\sum^m_{i=1}\frac{\partial \psi}{\partial x_i}(F_1,\cdots, F_m)D_tF_i.
\end{eqnarray}

 \item [(ii)] {\bf Integration by parts/duality formula} ([N], page 35)\\
Suppose $u(t)$ is $\cf_t-$adapted with $E[\int^T_0u^2(t)dt]<\infty$ and let $F\in \md^{(B)}_{1,2}$.
Then \begin{eqnarray}
E[F\int^T_0u(t)dB(t)]=E[\int^T_0u(t)D_tFdt].
\end{eqnarray}
\end{enumerate}

\subsection{Malliavin calculus for $\tilde N(\cdot)$}
The construction of a stochastic derivative/Malliavin derivative in the pure jump martingale
case follows the same lines as in the Brownian motion case. In this case the corresponding
Wiener-It\^{o} chaos expansion theorem states that any $F\in L^2({\cf}_T,P)$ (where in this case
$\cf_t=\cf^{(\tilde{N})}_t$ is the $\sigma-$algebra generated by $\eta(s):=\int^s_0\int_{\mr_0}
z\tilde{N}(dr,dz);\; 0\leq s\leq t$) can be written as
\begin{eqnarray}\label{WIcexpan}
F=\sum^\infty_{n=0}I_n(f_n);\; f_n\in \hat{L^2}((\lambda\times \nu)^n),
\end{eqnarray}
where $\hat{L^2}((\lambda\times \nu)^n)$ is the space of functions $f_n(t_1,z_1,\ldots, t_n,z_n)$;
$t_i\in[0,T$], $z_i\in \mr_0$ such that $f_n\in L^2((\lambda\times \nu)^n)$ and $f_n$ is symmetric
with respect to the pairs of variables $(t_1,z_1),\ldots,(t_n,z_n).$

 It is important to note that in this case the $n-$times iterated integral $I_n(f_n)$
 is taken with respect to $\tilde{N}(dt,dz)$ and not with respect to $d\eta(t).$
Thus, we define
\begin{equation}
I_n(f_n) =n!\int^T_0\!\!\int_{\mr_0}\!\int^{t_n}_0\!\int_{\mr_0}\cdots\int^{t_2}_0\!\!\int_{\mr_0}
f_n(t_1,z_1,\cdots,t_n,z_n)\tilde{N}(dt_1,dz_1)\cdots\tilde{N}(dt_n,dz_n)
\end{equation}
for $f_n\in \hat{L^2}((\lambda\times \nu)^n).$

The It\^{o} isometry for stochastic integrals with respect to $\tilde{N}(dt,dz)$
then gives the following isometry for the chaos expansion:
\begin{eqnarray}
||F||^2_{L^2(P)}=\sum^\infty_{n=0}n!||f_n||^2_{L^2((\lambda\times \nu)^n)}.
\end{eqnarray}
As in the Brownian motion case we use the chaos expansion to define the Malliavin derivative.
Note that in this case there are two parameters $t,z,$ where $t$ represents time and $z\neq 0$
represents a generic jump size.

\begin{defi}%\underline{Definition 2.3}
{\rm
{\bf (Malliavin derivative $D_{t,z}$)} [BDL\O{}P], [DM\O{}P]
Let $\md^{(\tilde{N})}_{1,2}$ be the space of all $F\in L^2({\cf}_T,P)$
such that its chaos expansion (\ref{WIcexpan}) satisfies
\begin{eqnarray}
||F||^2_{\md^{(\tilde{N})}_{1,2}}:=\sum^\infty_{n=1}n n!||f_n||^2_{L^2((\lambda\times \nu)^2)}<\infty.
\end{eqnarray}
For $F\in \md^{(\tilde{N})}_{1,2}$, we define the Malliavin derivative of $F$ at $(t,z)$ (with respect to
$\tilde{N(\cdot)})$, $D_{t,z}F, $ by
\begin{eqnarray}
D_{t,z}F=\sum^\infty_{n=1}nI_{n-1}(f_n(\cdot,t,z)),
\end{eqnarray}
where $I_{n-1}(f_n(\cdot,t,z))$ means that we perform the $(n-1)-$times iterated
integral with respect to $\tilde{N}$ to the first $n-1$ variable pairs $(t_1,z_1),\cdots,(t_n,z_n),
$ keeping $(t_n,z_n)=(t,z)$ as a parameter.}
\end{defi}

In this case we get the isometry.
\begin{eqnarray}
E[\int^T_0\int_{\mr_0}(D_{t,z}F)^2 \nu (dz)dt]=\sum^\infty_{n=0}n n!||f_n||^2_{L^2((\lambda\times \nu)^n)}=||F||^2_{\md_{1,2}^{(\tilde{N})}}.
\end{eqnarray}
(Compare with (\ref{isometry}).)
%\underline{Example 2.4}

\begin{ex}
\rm
If $F=\int^T_0\int_{\mr_0}f(t,z)\tilde{N}(dt,dz)$
for some deterministic $f(t,z)\in L^2(\lambda\times \nu)$, then
$$
D_{t,z}F=f(t,z) \mbox{ for } a. a. \, (t,z).
$$
More  generally, if  $\psi(s,\zeta)$ is Skorohod integrable with respect to
$\tilde N(\delta s, d\zeta)$, $\psi(s,\zeta)\in \md_{1,2}^{(\tilde N)}$ for
$a.a.\,s,\zeta$ and $D_{t,z}\psi(s, \xi)$ is Skorohod integrable for $a.a.\,(t,z)$, then
\begin{equation}
D_{t,z}(\int^T_0\!\int_\mr\psi(s,\zeta)\tilde{N}(\delta s,d\zeta))=\int^T_0\int_\mr D_{t,z}\psi(s,\zeta)\tilde{N}(\delta s,d\zeta)+u(t,z)\;
\mbox{ for } a.a.\, t,z,
\end{equation}
%for a.a. t,z,\\
where $\int^T_0\int_\mr\psi(s,\zeta)\tilde{N}(\delta s,d\zeta)$ denotes the
{\em Skorohod integral\/} of $\psi$ with respect to $\tilde{N}(\cdot,\cdot).$ (See [DM{\O}P] for a definition of such Skorohod integrals
and for more details.)
\end{ex}

The properties of $D_{t,z}$ corresponding to the properties (2.8) and (2.9) of $D_t$ are the following:
\begin{itemize}
\item [(i)] {\bf Chain rule([I], [DM{\O}P])} \quad
Suppose $F_1,\cdots, F_m\in \md^{(\tilde{N})}_{1,2}$ and that $\phi:\mr^m\rightarrow \mr$
is continuous and bounded. Then $\phi(F_1,\cdots,F_m)\in \md^{(\tilde{N})}_{1,2}$ and
\begin{equation}
D_{t,z}\phi(F_1,\cdots,F_m)=\phi(F_1+D_{t,z}F_1,\ldots, F_m+D_{t,z}F_m)-\phi(F_1,\ldots,F_m).
\end{equation}
\item [(ii)] {\bf Integration by parts/duality formula [DM{\O}P]}
\quad
Suppose $\Psi(t,z)$ is ${\cf}_t$-adapted and
$E[\int^T_0\int_{\mr_0}\psi^2(t,z)\nu(dz)dt]<\infty$
and let $F\in \md_{1,2}^{(\tilde{N})}$. Then
\begin{eqnarray}
E\Big[F\int^T_0\int_{\mr_0}\Psi(t,z)\tilde{N}(dt,dz)\Big]=
E\Big[\int^T_0\int_{\mr_0}\Psi(t,z)D_{t,z}F \nu(dz) dt\Big].
\end{eqnarray}
\end{itemize}

\section{The stochastic  maximum principle} % sec. 3

We now return to Problem \ref{prob1} given in the introduction. We make the following assumptions:
\begin{assump}
{\rm
\hspace{1cm}
\begin{itemize}
 \item [(3.1)] The functions  $b:[0,T]\times\mr\times U\times \Omega\rightarrow \mr$,
$\sigma:[0,T]\times\mr\times U\times \Omega \rightarrow \mr$,
$f:[0,T]\times \mr\times U\times \Omega \rightarrow \mr$ and
$g:\mr\times \Omega\rightarrow \mr$ are all continuously
differentiable $(C^1)$ with respect to $x\in \mr$ and $u\in U$ for
each $t \in [0,T]$ and a.a. $\omega\in\Omega$.
 \item [(3.2)] For all $t,r\in (0,T),\,t\le r$, and all bounded $\mE_t-$measurable
random variables $\alpha=\alpha(\omega)$ the control
$$\beta_\alpha(s)=\alpha(\omega)\chi_{[t,r]}(s);\quad s\in[0,T]$$
belongs to $\ma_\mE$.
 \item [ (3.3)] For all $u,\beta \in\ma_\mE$ with $\beta$ bounded, there exists $\delta>0$
such that
$$u+y\beta\in\ma_\mE\quad \mbox{ for all } y\in(-\delta,\delta)$$
and such that the family
\begin{align*}
&\Big\{ \frac{\partial f}{\partial x}(t,X^{u+y\beta}(t),u(t)+y\beta(t))
              \frac{d}{dy}X^{u+y\beta}(t) \\
&\qquad +\frac{\partial f}{\partial u}(t,X^{u+y\beta}(t),u(t)
        +y\beta(t))\beta(t)\Big\}_{y\in(-\delta,\delta)}
\end{align*}
is $\lambda\times P$-uniformly integrable and the family
$$
\Big\{
g'(X^{u+y\beta}(T))\frac{d}{dy}X^{u+y\beta}(T)\Big\}_{y\in(-\delta,\delta)}
$$
is $P$-uniformly integrable.
 \item [(3.4)] For all $u,\beta\in\ma_\mE$ with $\beta$ bounded the process
$Y(t)=Y^{(\beta)}(t)=\frac{d}{dy}X^{(u+y\beta)}(t)|_{y=0}$ exists and satisfies
the equation
\begin{align*}
&dY(t) =Y(t^-)\Big[\frac{\partial b}{\partial x}(t,X(t),u(t))dt
          +\frac{\partial \sigma}{\partial x}(t,X(t),u(t))dB(t) \\
&\qquad +\int_{\mr_0}\frac{\partial\theta}{\partial x}(t,X(t^-),u(t^-),z)\tilde{N}(dt,dz)\Big] \\
&\qquad +\beta(t^-)\Big[\frac{\partial b}{\partial u}(t,X(t),u(t))dt
       +\frac{\partial \sigma}{\partial u}(t,X(t),u(t))dB(t) \\
&\qquad +\int_{\mr_0}\frac{\partial \theta}{\partial u}(t,X(t^-),u(t^-),z)\tilde{N}(dt,dz)\Big];\\
&Y(0) =0.
\end{align*}

\item [(3.5)] For all $u\in\ma_\mE$, the following processes
\setcounter{equation}{5}
\begin{align}
K(t) &:=g'(X(T))+\int^T_t\frac{\partial f}{\partial x}(s,X(s),u(s))ds,\notag\\
D_tK(t) &:=D_tg'(X(T))+\int_t^TD_t\frac{\partial f}{\partial x}(s,X(s),u(s))ds, \notag\\
D_{t,z}K(t) &:=D_{t,z}g'(X(T))+\int^T_tD_{t,z}\frac{\partial f}{\partial x}(s,X(s),u(s))ds, \notag\\
H_0(s,x,u) &:=K(s)b(s,x,u)+D_sK(s)\sigma(s,x,u) \notag\\
&\qquad +\int_{\mr_0}D_{s,z}K(s)\theta(s,x,u,z)\nu(dz), \notag\\
G(t,s) &:= \exp \Big( \int_{t}^s \Big\{ \frac{\partial b}{\partial x}(r,X(r),u(r),\omega)
-\tfrac{1}{2}\Big( \frac{\partial\sigma}{\partial x}\Big)^2(r,X(r),u(r),\omega)\Big\} dr \notag\\
&\qquad +\int_{t}^s \frac{\partial\sigma}{\partial x}(r,X(r),u(r),\omega)dB(r) \notag\\
&\qquad +\int_{t}^s \int_{\mr_0}\Big\{ \ln \Big(
1+\frac{\partial\theta}{\partial x}(r,X(r),u(r),z,\omega)\Big)
     -\frac{\partial\theta}{\partial x}(r,X(r),u(r),z,\omega)\Big\} \nu(dz)dr \notag\\
&\qquad +\int_{t}^s \int_{\mr_0}\ln \Big(
1+\frac{\partial\theta}{\partial x}(r,X(r^-),u(r^-),z,\omega)\Big)
          \tilde{N}(dr,dz) \Big)\,,\notag\\
p(t)&:=K(t)+\int_t^T\frac{\partial H_0}{\partial x}(s,X(s),u(s))G(t,s)ds, \\
q(t)&:=D_tp(t)\,,\qquad \hbox{and} \\
r(t,z)&:=D_{t,z}p(t)
\end{align}
\end{itemize}
\hspace*{0cm}
all exist for $0\leq t\leq s\leq T,\;\, z\in\mr_0$.}
\end{assump}

We now define the {\em Hamiltonian\/} for this general problem:

\begin{defi}[The general stochastic Hamiltonian]%Definition 3.2
The general stochastic Hamiltonian is the process
$$
H(t,x,u,\omega):[0,T]\times\mr\times U\times\Omega\rightarrow\mr
$$
defined by
\begin{align}\nonumber
H(t,x,u,\omega)&=f(t,x,u,\omega)+p(t)b(t,x,u,\omega)+q(t)\sigma(t,x,u,\omega)\\
&+\int_{\mr_0}r(t,z)\theta(t,x,u,z,\omega)\nu(dz).
\end{align}
\end{defi}

\begin{remark}%Remark 3.3
\rm
In the classical Markovian case, the Hamiltonian
$H_1:[0,T]\times\mr\times U\times\mr\times\mr\times \cR \rightarrow \mr$ is defined by
\begin{equation}
H_1(t,x,u,p,q,r)=f(t,x,u)+p b(t,x,u)+q\sigma(t,x,u)+\int_{\mr_0}r(z)\theta(t,x,u,z)\nu(dz),
\end{equation}
where $\cR$ is the set of functions $r:\mr_0\rightarrow\mr$; see [F\O S]. Thus the relation between $H_1$ and $H$ is that:
\begin{equation}
H(t,x,u,\omega)=H_1(t,x,u,p(t),q(t),r(t,\cdot))
\end{equation}
where $p(t),q(t)$ and $r(t,z)$ are given by (3.6)--(3.8).
\end{remark}

We can now formulate our  stochastic maximum principle:

\begin{theo}[Maximum Principle]\label{GPIMPI}%Theorem 3.4
\hfill\break
{\bf (i)} \
Suppose $\hat{u}\in\ma_\mE$ is a critical point for $J(u)$, in the sense that
\begin{equation}\label{critpt}
\frac{d}{dy}J(\hat{u}+y\beta)|_{y=0}=0\quad \mbox{ for all bounded } \beta\in\ma_\mE.
\end{equation}
Then
\begin{equation}\label{conclu}
E[\frac{\partial\hat{H}}{\partial
u}(t,\hat{X}(t),\hat{u}(t))|\mE_t]=0 \mbox{ for  a.a. t,
$\omega$,}
\end{equation}
where
\begin{align*}
&\hat{X}(t)=X^{(\hat{u})}(t)\,, \\
&\hat{H}(t,\hat{X}(t),u)=f(t,\hat{X}(t),u)+\hat{p}(t)b
         (t,\hat{X}(t),u)+\hat{q}(t)\sigma(t,\hat{X}(t),u) \\
&\qquad +\int_{\mr_0}\hat{r}(t,z)\theta(t,\hat{X}(t),u,z)\nu(dz)\,,\\
\end{align*}
with
\begin{align*}
&\hat{p}(t)=\hat{K}(t)+\int_t^T \frac{\partial H_0}{\partial
x}(s,\hat{X}(s),\hat{u}(s))\hat{G}(t,s)ds\,,\\
\end{align*}
and
\begin{align*}
& \hat{G}(t,s) = \exp \Big( \int_{t}^s \Big\{ \frac{\partial
b}{\partial x}(r,\hat{X}(r),u(r),\omega)
-\tfrac{1}{2}\Big( \frac{\partial\sigma}{\partial x}\Big)^2(r,\hat{X}(r),u(r),\omega)\Big\} dr \notag\\
&\qquad +\int_{t}^s \frac{\partial\sigma}{\partial x}(r,\hat{X}(r),u(r),\omega)dB(r) \notag\\
&\qquad +\int_{t}^s \int_{\mr_0}\Big\{ \ln \Big(
1+\frac{\partial\theta}{\partial
x}(r,\hat{X}(r),u(r),z,\omega)\Big)
     -\frac{\partial\theta}{\partial x}(r,\hat{X}(r),u(r),z,\omega)\Big\} \nu(dz)dr \notag\\
&\qquad +\int_{t}^s \int_{\mr_0}\ln \Big(
1+\frac{\partial\theta}{\partial
x}(r,\hat{X}(r^-),u(r^-),z,\omega)\Big) \tilde{N}(dr,dz) \Big)\,,
\\
&\hat{K}(t)=K^{(\hat{u})}(t)=g'(\hat{X}(T))+\int^T_t\frac{\partial
f}{\partial x}(s,\hat{X}(s),\hat{u}(s))ds.
\end{align*}
{\bf (ii)}æ\
Conversely, suppose there exists $\hat{u}\in\ma_\mE$ such that (3.13) holds. Then $\hat{u}$ satisfies (3.12).
\end{theo}

\noindent
{\sc Proof.} \\
{\bf (i):} \
 Suppose $\hat{u}\in\ma_\mE$ is a critical point for $J(u)$. Choose an arbitrary $\beta\in\ma_\mE$ bounded and let $\delta >0$ be as in
(3.3) of Assumption 3.1.

For simplicity of notation  we write $\hat{u}=u,\hat{X}=X \mbox{ and } \hat{Y}=Y$ in the following. By (3.3) we have
\begin{equation}
\begin{array}{ll}
0&=\frac{d}{dy}J(u+y\beta)|_{y=0} \\
&=E[\int^T_0\{\frac{\partial f}{\partial x}(t,X(t),u(t))Y(t)+\frac{\partial f}{\partial u}(t,X(t),u(t))\beta(t)\}dt+g'(X(T))Y(T)],
\end{array}
\end{equation}
where
\begin{equation}
\begin{array}{ll}
Y(t)&=Y^{(\beta)}(t)=\frac{d}{dy}X^{(u+y\beta)}(t)|_{y=0}\\
&=\int^t_0\{\frac{\partial b}{\partial x}(s,X(s),u(s))Y(s)+\frac{\partial b}{\partial u}(s,X(s),u(s))\beta(s)\}ds\\
&+\int^t_0\{\frac{\partial \sigma}{\partial x}(s,X(s),u(s))Y(s)+\frac{\partial \sigma}{\partial u}(s,X(s),u(s))\beta(s)\}dB(s)\\
&+\int^t_0\int_{\mr_0}\{\frac{\partial \theta}{\partial x}(s,X(s),u(s),z)Y(s)+\frac{\partial\theta}{\partial u}(s,X(s),u(s),z)\beta(s)\}\tilde{N}(ds,dz).
\end{array}
\end{equation}

If we use the short hand notation
$$\frac{\partial f}{\partial x}(t,X(t),u(t))=\frac{\partial f}{\partial x}{(t)},\;\; \frac{\partial f}{\partial u}(t,X(t),u(t))=\frac{\partial f}{\partial u}(t)$$
and similarly for $\frac{\partial b}{\partial x},\frac{\partial b}{\partial u},\frac{\partial \sigma}{\partial x},\frac{\partial \sigma}{\partial u},\frac{\partial \theta}{\partial x}$, and $\frac{\partial \theta}{\partial u}$,
we can write
\begin{equation}
\begin{array}{ll}
dY(t)&=\{\frac{\partial b}{\partial x}(t)Y(t)+\frac{\partial b}{\partial u}(t)\beta(t)\}dt
      +\{\frac{\partial \sigma}{\partial x}(t)Y(t)+\frac{\partial \sigma}{\partial u}(t)
      \beta(t)\}dB(t)\\
&+\int_{\mr_0}\{\frac{\partial \theta}{\partial x}(t)Y(t)+\frac{\partial \theta}{\partial u}(t)
         \beta(t)\}\tilde{N}(dt,dz);\\
Y(0)&=0.
\end{array}
\end{equation}
By the duality formulas (2.9) and (2.18), we get
\begin{align}
%\begin{array}{ll}
E[g' & (X(T))Y(T)] \nonumber \\
&={\textstyle E\Big[g'(X(T))\Big(\int^T_0\{\frac{\partial b}{\partial x}(t)Y(t)
       +\frac{\partial b}{\partial u}(t)\beta(t)\}dt} \nonumber \\
&+{\textstyle \int^T_0\{\frac{\partial \sigma}{\partial x}(t)Y(t)
             +\frac{\partial\sigma}{\partial u}(t)\beta(t)\}dB(t)} \nonumber  \\
&+{\textstyle \int^T_0\int_{\mr_0}\{\frac{\partial\theta}{\partial x}(t)Y(t)+
    \frac{\partial \theta}{\partial u}(t)\beta(t)\}\tilde{N}(dt,dz)\Big)\Big]} \nonumber  \\
&={\textstyle E\Big[\int^T_0\{g'(X(T))\frac{\partial b}{\partial x}(t)Y(t)
         +g'(X(T))\frac{\partial b}{\partial u}(t)\beta(t)} \nonumber  \\
&+{\textstyle D_t(g'(X(T)))\frac{\partial \sigma}{\partial x}(t)Y(t)+D_t(g'(X(T)))
       \frac{\partial \sigma}{\partial u}(t)\beta(t)} \nonumber  \\
&+{\textstyle \int_{\mr_0}[D_{t,z}(g'(X(T)))\frac{\partial \theta}{\partial x}(t)Y(t)+D_{t,z}(g'(X(T)))
           \frac{\partial\theta}{\partial u}(t)\beta(t)]\nu(dz)\}dt\Big].}
%\end{array}
\end{align}
Similarly we have, using the Fubini theorem,
\begin{align*}
&{\textstyle E[\int^T_0\frac{\partial f}{\partial x}(t)Y(t)dt]} \\
&\quad={\textstyle E[\int^T_0\frac{\partial f}{\partial x}(t)(\int^t_0\{\frac{\partial b}{\partial x}(s)Y(s)
               +\frac{\partial b}{\partial u}(s)\beta(s)\}ds} \\
&\qquad {\textstyle +\int^t_0\{\frac{\partial \sigma}{\partial x}(s)Y(s)
                +\frac{\partial \sigma}{\partial u}(s)\beta(s)\}dB(s)} \\
&\qquad {\textstyle +\int^t_0\int_{\mr_0}\{\frac{\partial \theta}{\partial x}(s)Y(s)
              +\frac{\partial\theta}{\partial u}(s)\beta(s)\}\tilde{N}(ds,dz))dt]} \\
&\quad {\textstyle =E[\int^T_0(\int^t_0\{\frac{\partial f}{\partial x}(t)[\frac{\partial b}{\partial x}(s)Y(s)
                 +\frac{\partial b}{\partial u}(s)\beta(s)]} \\
&\qquad {\textstyle +D_s(\frac{\partial f}{\partial x}(t))[\frac{\partial \sigma}{\partial x}(s)Y(s)
               +\frac{\partial \sigma}{\partial u}(s)\beta(s)]} \\
&\qquad {\textstyle +\int_{\mr_0}D_{s,z}(\frac{\partial f}{\partial x}(t))[
            \frac{\partial\theta}{\partial x}(s)Y(s)
        +\frac{\partial\theta}{\partial u}(s)\beta(s)]\nu(dz)\}ds)dt} \\
&\quad {\textstyle =E[\int_0^T\{ ( \int_s^T \frac{\partial f}{\partial x}(t)dt)[
            \frac{\partial b}{\partial x}(s)Y(s)
        +\frac{\partial b}{\partial u}(s)\beta(s)]} \\
&\qquad {\textstyle +(\int_s^T D_s\frac{\partial f}{\partial x}(t)dt)
        [\frac{\partial\sigma}{\partial x}(s)Y(s)+\frac{\partial\sigma}{\partial u}(s)\beta(s)]} \\
&\qquad {\textstyle +\int_{\mr_0}(\int^T_sD_{s,z}\frac{\partial f}{\partial x}(t)dt)
        [\frac{\partial \theta}{\partial x}(s)Y(s)
              +\frac{\partial \theta}{\partial u}(s)\beta(s)]\nu(dz)\}ds].}
\end{align*}
Changing the notation $s\rightarrow t$, this becomes
\begin{equation}
\begin{array}{ll}
E[\int^T_0\frac{\partial f}{\partial x}(t)Y(t)dt]=&E[\int^T_0\{(\int^T_t\frac{\partial f}{\partial x}(s)ds)[\frac{\partial b}{\partial x}(t)Y(t)
               +\frac{\partial b}{\partial u}(t)\beta(t)]\\
&+(\int^T_tD_t\frac{\partial f}{\partial x}(s)ds)[\frac{\partial \sigma}{\partial x}(t)Y(t)
                +\frac{\partial \sigma}{\partial u}(t)\beta(t)]\\
&+\int_{\mr_0}(\int^T_tD_{t,z}\frac{\partial f}{\partial x}(s)ds)[\frac{\partial\theta}{\partial x}(t)Y(t)
               +\frac{\partial \theta}{\partial u}(t)\beta(t)]\nu(dz)\}dt].
\end{array}
\end{equation}
Recall
\begin{equation}\label{end}
K(t):=g'(X(T))+\int^T_t\frac{\partial f}{\partial x}(s)ds.
\end{equation}
By combining (3.17)--(3.19), we get
\begin{equation}\label{to0}
\begin{array}{ll}
&E[\int^T_0\{K(t)(\frac{\partial b}{\partial x}(t)Y(t)+\frac{\partial b}{\partial u}(t)\beta(t))\\
+&D_tK(t)(\frac{\partial \sigma}{\partial x}(t)Y(t)+\frac{\partial \sigma}{\partial u}(t)\beta(t))\\
+&\int_{\mr_0}D_{t,z}K(t)(\frac{\partial \theta}{\partial x}(t)Y(t)+\frac{\partial \theta}{\partial u}(t)\beta(t))\nu(dz)+\frac{\partial f}{\partial u}(t)\beta(t)\}dt]=0.
\end{array}
\end{equation}

Now apply the above to $\beta=\beta_\alpha\in\ma_\mE $ of the form $\beta_\alpha(s)=\alpha\chi_{[t,t+h]}(s),$ for some $t,h\in(0,T)$, $t+h\leq T$,
where $\alpha=\alpha(\omega)$ is bounded and $\mE_t$-measurable.
Then $Y^{(\beta_\alpha)}(s)=0$ for $0\leq s\leq t$ and hence (3.20) becomes
\begin{equation}
A_1+A_2=0,
\end{equation}
where
\begin{align*}
A_1&=E\left[\int^T_t\{K(s)\frac{\partial b}{\partial x}(s)
       +D_sK(s)\frac{\partial \sigma}{\partial x}(s)+\int_{\mr_0}D_{s,z}K(s)
       \frac{\partial \theta}{\partial x}(s)\nu(dz)\}Y^{(\beta_\alpha)}(s)ds\right],\\
A_2&=E\left[(\int^{t+h}_t\!\{K(s)\frac{\partial b}{\partial u}(s)
      +D_sK(s)\frac{\partial \sigma}{\partial u}(s)
      +\int_{\mr_0}\!D_{s,z}K(s)\frac{\partial\theta}{\partial u}(s)\nu(dz)
        +\frac{\partial f}{\partial u}(s)\}ds)\alpha\right].
\end{align*}
Note that, by (3.16), with $Y(s)=Y^{(\beta_\alpha)}(s)$ and $s\geq
t+h$, the process $Y(s)$ follows the following dynamics
\begin{equation}
dY(s)=Y(s^-)\Big\{ \frac{\partial b}{\partial
x}(s)ds+\frac{\partial\sigma}{\partial x}(s)dB(s) +\int_{\mr_0}
\frac{\partial\theta}{\partial x}(s)\tilde{N}(ds,dz)\Big\}\,,
\end{equation}
for $s\ge t+h$ with initial condition $Y(t+h)$ in time $t+h$. This
equation can be solved explicitly and we get
\begin{equation} % eq3.22
Y(s)=Y(t+h)G(t+h,s);\qquad s\geq t+h,
\end{equation}
where, in general, for $s\geq t$,
\begin{align*} % eq 3.23
&G(t,s) = \exp \Big( \int_{t}^s \Big\{ \frac{\partial b}{\partial x}(r)-\tfrac{1}{2}
        \Big( \frac{\partial\sigma}{\partial x}\Big)^2(r)\Big\}dr
          +\int_{t}^s \frac{\partial\sigma}{\partial x}(r)dB(r) \nonumber \\
&\qquad +\int_{t}^s\int_{\mr_0}\ln \Big( 1+\frac{\partial\theta}{\partial x}(r)\Big)
            \tilde{N}(dr,dz) \nonumber \\
&\qquad + \int_{t}^s \int_{\mr_0}æ\Big\{ \ln \Big(
       1+\frac{\partial\theta}{\partial x}(r)\Big)
         -\frac{\partial\theta}{\partial x}(r)\Big\} \nu(dz)dr\Big).
\end{align*}
That $Y(s)$ indeed is the solution of (3.22) can be verified
by applying the It\^o formula to $Y(s)$ given in (3.23). Note that
$G(t,s)$ does not depend on $h$, but $Y(s)$ does.

Put
\begin{equation}  % eq3.24
H_0(s,x,u)=K(s)b(s,x,u)+D_sK(s)\sigma(s,x,u)
+\int_{\mr_0} D_{s,z}K(s)\theta(s,x,u,z)\nu(dz),
\end{equation}
and $H_0(s)=H_0^{(u)}(s)=H_0(s,X(s),u(s))$. Then
$$
A_1=E\Big[ \int_t^T \frac{\partial H_0}{\partial x}(s)Y(s)ds\Big].
$$
Differentiating with respect to $h$ at $h=0$ we get
\begin{equation} % eq3.25
\frac{d}{dh} A_1\big|_{h=0} = \frac{d}{dh}E \Big[
\int_t^{t+h} \frac{\partial H_0}{\partial x}(s)Y(s)ds\Big]_{h=0} + \frac{d}{dh}E
\Big[ \int_{t+h}^T \frac{\partial H_0}{\partial x}(s)Y(s)ds\Big]_{h=0}.
\end{equation}
Since $Y(t)=0$ and since $\partial H_0/ \partial x (s)$ is c\`adl\`ag we see that
\begin{equation}
\frac{d}{dh}E \Big[ \int_t^{t+h} \frac{\partial H_0}{\partial x}(s)Y(s)ds\Big]_{h=0}=0.
\end{equation}
Therefore, using (3.23) and that $Y(t)=0$,
\begin{align}
\frac{d}{dh}A_1\big|_{h=0} &= \frac{d}{dh}E\Big[ \int_{t+h}^T
       \frac{\partial H_0}{\partial x}(s)Y(t+h)G(t+h,s)ds\Big]_{h=0} \nonumber \\
&= \int_t^T \frac{d}{dh}E\Big[ \frac{\partial H_0}{\partial x}(s)
                Y(t+h)G(t+h,s)\Big]_{h=0} ds \nonumber \\
&= \int_t^T \frac{d}{dh}E \Big[ \frac{\partial H_0}{\partial x}(s)G(t,s)Y(t+h)\Big]_{h=0}ds.
\end{align}
By (3.16)
\begin{align}
Y(t+h) &= \alpha \int_t^{t+h} \Big\{ \frac{\partial b}{\partial u}(r)dr
+\frac{\partial\sigma}{\partial u}(r)dB(r)
+\int_{\mr_0} \frac{\partial\theta}{\partial u}(r)\tilde{N}(dr,dz)\Big\} \nonumber \\
&\qquad +\int_t^{t+h}Y(r^-)\Big\{ \frac{\partial b}{\partial x}(r)dr
+\frac{\partial\sigma}{\partial x}(r)dB(r)
+\int_{\mr_0} \frac{\partial\theta}{\partial x}(r)\tilde{N}(dr,dz)\Big\}.
\end{align}
Therefore, by (3.27) and (3.28),
\begin{equation}
\frac{d}{dh}A_1\big|_{h=0}=\Lambda_1+\Lambda_2,
\end{equation}
where
\begin{align}
\Lambda_1 &=\int_t^T \frac{d}{dh}E \Big[ \frac{\partial H_0}{\partial x}(s)G(t,s)\alpha
        \int_t^{t+h} \Big\{ \frac{\partial b}{\partial u}(r)dr
        +\frac{\partial\sigma}{\partial u}(r)dB(r) \nonumber \\
&\qquad +\int_{\mr_0} \frac{\partial\theta}{\partial u}(r)\tilde{N}(dr,dz)\Big\}\Big]_{h=0}ds
\end{align}
and
\begin{align}
\Lambda_2 &=\int_t^T \frac{d}{dh}E \Big[ \frac{\partial H_0}{\partial x}(s)G(t,s)
             \int_t^{t+h} Y(r^-)\Big\{ \frac{\partial b}{\partial x}(r)dr
        +\frac{\partial\sigma}{\partial x}(r)dB(r) \nonumber \\
&\qquad +\int_{\mr_0} \frac{\partial\theta}{\partial x}(r)\tilde{N}(dr,dz)\Big\}\Big]_{h=0}ds.
\end{align}
By the duality formulae (2.9), (2.18) we have
\begin{align}
\Lambda_1 &=\int_t^T \frac{d}{dh}E \Big[ \alpha\int_t^{t+h} \Big\{
       \frac{\partial b}{\partial u}(r)F(t,s)
        +\frac{\partial\sigma}{\partial u}(r)D_rF(t,s) \nonumber \\
&\qquad +\int_{\mr_0} \frac{\partial\theta}{\partial u}(r)D_{r,z}F(t,s)
             \nu(dz)\Big\}dr\Big]_{h=0}ds \nonumber \\
&=\int_t^T E \Big[ \alphaæ\Big\{ \frac{\partial b}{\partial u}(t)F(t,s)
               +\frac{\partial\sigma}{\partial u}(t)D_tF(t,s)
       +\int_{\mr_0}\frac{\partial\theta}{\partial u}(t)D_{t,z}F(t,s)\nu(dz)\Big\}\Big]ds,
\end{align}
where we have put
\begin{equation}
F(t,s)=\frac{\partial H_0}{\partial x}(s)G(t,s).
\end{equation}
Since $Y(t)=0$  we see that
\begin{equation}
\Lambda_2=0.
\end{equation}
We conclude that
\begin{align}
&\frac{d}{dh}A_1\big|_{h=0} = \Lambda_1 \nonumber \\
&\quad= \int_t^T E\Big[ \alpha \Big\{ F(t,s)\frac{\partial b}{\partial u}(t)
+D_tF(t,s)\frac{\partial\sigma}{\partial u}(t)
+\int_{\mr_0} D_{t,z}F(t,s)\frac{\partial\theta}{\partial u}(t)\nu(dz)\Big\}\Big]ds.
\end{align}
Moreover, we see directly that
\[
\frac{d}{dh}A_2\big|_{h=0}
= E\Big[ \alpha \Big\{ K(t)\frac{\partial b}{\partial u}(t)
       +D_tK(t)\frac{\partial\sigma}{\partial u}(t)
         +\int_{\mr_0} D_{t,z}K(t)\frac{\partial\theta}{\partial u}(t)\nu(dz)
+\frac{\partial f}{\partial u}(t)\Big\}\Big].
\]
Therefore, differentiating (3.21) with respect to $h$ at $h=0$ gives the equation
\begin{align}
&E\Big[ \alpha \Big\{ \Big( K(t)+\int_t^T F(t,s)ds\Big) \frac{\partial b}{\partial u}(t)
       +D_t\Big( K(t)+\int_t^T F(t,s)ds\Big) \frac{\partial\sigma}{\partial u}(t) \nonumber \\
&\qquad +\int_{\mr_0} D_{t,z}\Big( K(t)+\int_t^T F(t,s)ds\Big)
      \frac{\partial\theta}{\partial u}(t)\nu(dz+\frac{\partial f}{\partial u}(t))\Big\}\Big]=0.
\end{align}
We can reformulate this as follows: If we define, as in (3.6),
\begin{equation}
p(t)=K(t)+\int_t^T F(t,s)ds=K(t)+\int_t^T \frac{\partial H_0}{\partial x}(s)G(t,s)ds,
\end{equation}
then (3.36)  can be written
\begin{align*}
E\Big[ &\frac{\partial}{\partial u}\Big\{ f(t,X(t),u)+p(t)b(t,X(t),u)+D_t p(t)\sigma(t,X(t),u)  \\
&+\int_{\mr_0} D_{t,z}p(t)\theta(t,X(t),u,z)\nu(dz)\Big\}_{u=u(t)}\alpha\Big]=0.
\end{align*}
Since this holds for all bounded $\mE_t$-measurable random variable $\alpha$, we conclude that
$$
E\Big[ \frac{\partial}{\partial u}H(t,X(t),u)_{u=u(t)}\mid \mE_t\Big]=0,
$$
which is (3.13). This completes the proof of (i).
\bigskip

\noindent
{\bf (ii):}æ\
Conversely, suppose (3.13) holds for some $\hat{u}\in\ma_\mE$. Then by reversing the above argument we get that (3.21) holds for all $\beta_\alpha\in\ma_\mE$ of the form
$$
\beta_\alpha(s,\omega)=\alpha(\omega)\chi_{(t,t+h]}(s)
$$
for some $t,h\in[0,T]$ with $t+h\leq T$ and some bounded $\mE_t$-measurable $\alpha$. Hence (3.21) holds for all linear combinations of such $\beta_\alpha$. Since all bounded $\beta\in\ma_\mE$ can be approximated pointwise boundedly in $(t,\omega)$ by such linear combinations, it follows that (3.21) holds for all bounded $\beta\in\ma_\mE$. Hence, by reversing the remaining part of the argument above, we conclude that (3.12) holds.

\hfill
\eof

\section{Applications}

In this section we illustrate the maximum principle by looking at some examples.

\begin{ex}[Optimal dividend/harvesting rate]
\hfill\break
\rm
Suppose the cash flow $X(t)=X^{(c)}(t)$ at time $t$ is given by
\begin{equation}
\begin{array}{l}
dX(t)=(b_0(t,\omega)+b_1(t,\omega)X(t)-c(t))dt+(\sigma_0(t,\omega)+\sigma_1(t,\omega)X(t))dB(t)\\
\hspace{5.7cm}+\int_{\R_0}(\theta_0(t,z,\omega)+\theta_1(t,z,\omega)X(t))\tilde N(dt,dz)\,; \\
X(0)=x\in \R,
\end{array}
\end{equation}
where
$$
\begin{array}{ll}
b_0(t)=b_0(t,\omega),\ b_1(t)=b_1(t,\omega):              & [0,T]\times \Omega \mapsto \R\\
\sigma_0(t)=\sigma_0(t,\omega),\ \sigma_1(t)=\sigma_1(t,\omega):     & [0,T]\times \Omega \mapsto \R\;\; \rm{and}\\
\theta_0(t,z)=\theta_0(t,z,\omega),\
\theta_1(t,z)=\theta_1(t,z,\omega):   & [0,T]\times\R \times\Omega
\mapsto \R
\end{array}$$
are given $\cf_t$-predictable processes.

Here $c(t)\ge 0$ is our control (the dividend/harvesting rate), assumed to belong to a family ${\cal A}_\mE$ of admissible controls, contained in the set of $\mE_t$-predictable controls.

Suppose the performance functional has the form
\begin{equation}
J(c)=E\left[\int_0^T\xi(s)U(c(s))ds+\zeta X^{(c)}(T)\right]
\end{equation}
where $U: [0,+\infty]\mapsto \R$ is a $C^1$ utility function, $\xi(s)=\xi(s,\omega)$ is an $\cf_t$-predictable process and $\zeta=\zeta(\omega)$ is an $\cf_T$-measurable random variable.

We want to find $\hat c\in {\cal A}_\mE$ such that
\begin{equation}\label{probex41}
 \sup_{c\in {\cal A}_\mE}J(c)=J(\hat c).
\end{equation}

Using the notation from the previous section, we note that in this
case we have, with $c=u$,
$$
f(t,x,c)=\xi(t)U(c)\text{\ \ \ \ and \ \ \ \ } g(x)=\zeta x.
$$
Hence
\begin{eqnarray*}
K(t)&=&\int_t^T\frac{\partial f}{\partial x}(s,X(s),c(s))ds+g'(X(T))=\zeta\,,\\
H_0(t,x,c)&=& \zeta (b_0(t)+b_1(t)x-c) + D_t\zeta(\sigma_0(t)+\sigma_1(t)x) + \int_{\mr_0}D_{t,z}\zeta(\theta_0(t,z)+\theta_1(t,z)x)\nu(dz)\,, \notag\\
G(t,s)&=&\exp\Big(\int_{t}^s\Big\{b_1(r)-\tfrac{1}{2}\sigma_1^2(r)\Big\}\,dr
    +\int_{t}^s \sigma_1(r)\,dB(r) \notag\\
&&+\int_{t}^s \int_{\mr_0}\Big\{ \ln \Big( 1+\theta_1(r,z)\Big)
     -\theta_1(r,z)\Big\} \nu(dz)dr+\int_{t}^s \int_{\mr_0}\ln \Big( 1+\theta_1(r,z)\Big)\tilde{N}(dr,dz)\Big)
     \,.
\end{eqnarray*}
Then
\begin{eqnarray}
\label{eq:p(t)} p(t)&=& \zeta + \int_{t}^T \left( \zeta b_1(r) +
D_r\zeta\,\sigma_1(r)+\int_{\mr_0}D_{r,z}\zeta\,\theta_1(r,z)\nu(dz)\right)G(t,r)\,dr\,,
\end{eqnarray}
and the Hamiltonian becomes
\begin{eqnarray}
H(t,x,c)&=&\xi(t)U(c)+p(t)(b_0(t)+b_1(t)x-c)+D_tp(t)\,(\sigma_0(t)+\sigma_1(t)x)\\
&&\hspace{4cm}+\int_{\mr_0}D_{t,z}p(t)\,(\theta_0(t,z)+\theta_1(t,z)x)\nu(dz)
\notag.
\end{eqnarray}
Hence, if $\hat c\in {\cal A}_\mE$ is optimal for the problem (\ref{probex41}), we have
\begin{eqnarray*}
0&=&E\left[\frac{\partial}{\partial c}H(t,\hat X(t), \hat c(t))|\mE_t\right]\\
&=&E\left[\{\xi(t)U'(\hat c(t))-p(t)\}|\mE_t\right]\\
&=&U'(\hat c(t))E\left[\xi(t)|\mE_t\right]-E[p(t)|\mE_t].
\end{eqnarray*}
\end{ex}

We have proved:

\begin{theo}
If there exists an optimal dividend/harvesting rate $\hat c(t)>0$ for problem (\ref{probex41}), then it satisfies the equation
\begin{equation}
U'(\hat c(t))E[\xi(t)|\mE_t]=E[p(t) |\mE_t]\,,
\end{equation}
where $p(t)$ is given by \eqref{eq:p(t)}.
\end{theo}

\begin{ex}[Optimal portfolio]
\hfill\break
\rm
Suppose we have a financial market with the following two investment possibilities:
\medskip

\noindent
(i)\
A risk free asset, where the unit price $S_0(t)$ at time $t$ is given by
\begin{equation} % 4.7
\begin{array}{l}
dS_0(t)=\rho_tS_0(t)dt;\qquad S_0(0)=1;\quad t\in[0,T]\,. \\
\end{array}
\end{equation}
(ii) \
A risky asset, where the unit price $S_1(t)$ at time $t$ is given by
\begin{equation}
\begin{array}{l} % 4.8
dS_1(t)=S_1(t^-)\Big[ \alpha_t dt+\beta_t dB(t)
+\int_{\mr_0}\zeta(t,z)\tilde{N}(dt,dz)\Big];\qquad t\in[0,T]\\
S_1(0)>0.
\end{array}
\end{equation}
\end{ex}
Here $\rho_t,\alpha_t,\beta_t$ and $\zeta(t,z)$ are bounded $\cf_t$-predictable processes;
$t\in[0,T]$, $z\in\mr_0$ and $T>0$ is a given constant. We also assume that
$$
\zeta(t,z)\geq -1\qquad \hbox{a.s. for a.a. $t,z$}
$$
and
$$
E\Big[ \int_0^T \int_{\mr_0} | \log(1+\zeta(t,z))|^2 \nu(dz)dt\Big]<\infty\,.
$$
A {\em portfolio\/} in this market is an $\ce_t$-predictable
process $u(t)$ representing the amount invested in the risky asset
at time $t$. When the portfolio $u(\cdot)$ is chosen, the
corresponding wealth process $X(t)=X^{(u)}(t)$ satisfies the
equation
\begin{align}
dX(t)=[  & \rho_t X(t)+(\alpha_t-\rho_t)u(t)]dt+\beta_t u(t)dB(t) \nonumber \\
&+\int_{\mr_0} \zeta(t,z)u(t^-)\tilde{N}(dt,dz);\qquad X(0)=x>0.
\end{align}
{\em The partial information optimal portfolio problem\/} is to find the portfolio $u\in \ca_\ce$
which maximizes
$$
J(u)=E[U(X^{(u)}(T),\omega)]
$$
where $U(x)=U(x,\omega):\mr\times \Omega\to\mr$ is a given $\cf_t$-measurable random variable for each $x$ and $x\to U(x,\omega)$ is a utility function for each $\omega$. We assume that $x\to U(x)$ is $C^1$ and $U'(x)$ is strictly decreasing. The set $\ca_\ce$ of admissible portfolios is contained in the set of $\ce_t$-adapted portfolios $u(t)$ such that (4.9) has a unique strong solution.

With the notation of the previous section we see that in this case we have
\begin{align*}
&f(t,x,u)=f(t,x,u,\omega)=0,\qquad g(x,\omega)=U(x,\omega), \\
&b(t,x,u)=\rho_t x+(\alpha_t-\rho_t)u,\qquad \sigma(t,x,u)=\beta_t u, \\
&\theta(t,x,u,z)=\zeta(t,z)u.
\end{align*}
Thus
\begin{equation} % 4.10
K(t)=U'(X(T))=K\,,
\end{equation}
and
\begin{align*}
H_0(t,x,u)=K & (\rho_t x+(\alpha_t-\rho_t)u)+D_tK\beta_t u  \\
&+\int_{\mr_0} D_{t,z}K\zeta(t,z)u\tilde{N}(dt,dz)\,, \nonumber
\end{align*}
and
\begin{equation} % eq4.12
G(t,s)=\exp\left(\int_t^s \rho_r dr\right).
\end{equation}
Thus
\begin{equation} % 4.13
p(t)=U'(X(T))+\int_t^T K\rho_s \exp\left(\int_t^s \rho_r dr\right)
ds\,,
\end{equation}
and the Hamiltonian becomes
\begin{align} % 4.14
H(t,x,u)=p(t) & [\rho_t x+(\alpha_t-\rho_t)u]+D_t p(t)\beta_t u  \\
&\hspace{2cm}+\int_{\mr_0} D_{t,z}p(t)\zeta(t,z)u\,\nu(dz).
\nonumber
\end{align}
By the maximum principle, we now get the following condition for
an optimal control.
\begin{theo}
If $\hat{u}(t)$ is an optimal control with corresponding
$\hat{X}(t),\hat{p}(t)$ then
$$
E\Big[ \frac{d}{du}H(t,\hat{X}(t),u)_{u=\hat{u}(t)}\mid \ce_t\Big]=0
$$
i.e.
\begin{equation}
E\Big[ \hat{p}(t)(\alpha_t-\rho_t)+\beta_t
D_t\hat{p}(t)+\int_{\mr_0}
D_{t,z}\hat{p}(t)\zeta(t,z)\,\nu(dz)\mid \ce_t\Big]=0.
\end{equation}
\end{theo}
Equation (4.14) is an interesting new type of equation. We could
call it a {\em Malliavin differential type equation\/} in the
unknown process $\hat{p}(t)$. Note that if we can find
$\hat{p}(t)$, then we also know $U'(\hat{X}(T))$ and hence
$\hat{X}(T)$. In particular, we see that the
optimal final wealth fulfills the following Malliavin differential
type equation in the unknown random variable $U'(\hat{X}(T))$
\begin{equation*}
E\Big[ U'(\hat{X}(T))(\alpha_T-\rho_T)+\beta_T
D_TU'(\hat{X}(T))+\int_{\mr_0}
D_{T,z}U'(\hat{X}(T))\zeta(T,z)\,\nu(dz)\mid \ce_T\Big]=0.
\end{equation*}

In this paper we will not discuss general solution methods of
this type of Malliavin differential equations, but leave this issue for future research.
Instead we complete by considering a solution in the special case when
\begin{equation}
\nu=\rho_t=0,\quad |\beta_t|\geq\delta>0\quad \hbox{and}\quad \ce_t=\cf_t;\;\; 0\leq t\leq T,
\end{equation}
where $\delta>0$ is a given constant. Then (4.14) simplifies to
\begin{equation}
\alpha_t E[K| \cf_t]+\beta_t E[D_tK| \cf_t]=0\,.
\end{equation}
By the Clark-Ocone theorem we have
$$
K=E[K]+\int_0^T E[D_tK| \cf_t]dB(t),
$$
which implies that
\begin{equation}
E[K| \cf_t]=E[K]+\int_0^t E[D_sK| \cf_s]dB(s).
\end{equation}
Define
\begin{equation}
M_t:=E[K| \cf_t]=E[U'(\hat{X}(T))| \cf_t].
\end{equation}
Then by substituting (4.16) into (4.17) we get
$$
M_t=E[K]-\int_0^t \frac{\alpha_s}{\beta_s}M_s dB(s)
$$
or
$$
dM_t=-\frac{\alpha_t}{\beta_t}M_t dB_t
$$
which has the solution
\begin{equation}
M_t=E[U'(\hat{X}(T))]\exp \Big( -\int_0^t \frac{\alpha_s}{\beta_s}dB(s)
-\tfrac{1}{2}\int_0^t\Big( \frac{\alpha_s}{\beta_s}\Big)^2 ds\Big).
\end{equation}
This determines $U'(\hat{X}(T))=M_T=K$ modulo the constant $E[U'(\hat{X}(T))]=M_0$.

Given $K$ the corresponding optimal portfolio $\hat{u}$ is given as the solution of the backward stochastic differential equation
\begin{equation}
\begin{cases}
d\hat{X}(t)=\alpha_t\hat{u}(t)dt+\beta_t\hat{u}(t)dB(t);\quad t<T \\
\hat{X}(T)=(U')^{-1}(K)\end{cases}
\end{equation}
This equation can be written
\begin{equation}
\begin{cases}
d\hat{X}(t)=\beta_t\hat{u}(t)d\tilde{B}(t);\quad t<T \\
\hat{X}(T)=(U')^{-1}(K),\end{cases}
\end{equation}
where
\begin{equation} % 4.23
d\tilde{B}(t)=\frac{\alpha_t}{\beta_t}dt+dB(t),
\end{equation}
which is a Brownian motion with respect to the probability measure $Q$ defined by
\begin{equation}
dQ=N_TdP\quad \hbox{on $\cf_T$\,,}
\end{equation}
where
\begin{equation*}
N_t=\exp \Big( -\int_0^t \frac{\alpha_s}{\beta_s}dB(s)
-\tfrac{1}{2}\int_0^t\Big( \frac{\alpha_s}{\beta_s}\Big)^2 ds\Big).
\end{equation*}
By the Clark-Ocone theorem under change of measure [KO] we have
\begin{equation} % 4.24
\hat{X}(T)=E_Q[\hat{X}(T)] + \int_0^T E_Q \Big[ \Big( D_t\hat{X}(T)
-\hat{X}(T)\int_t^T D_t \Big( \frac{\alpha_s}{\beta_s}\Big) d\tilde{B}(s)\Big) \mid \cf_t \Big]
d\tilde{B}(t).
\end{equation}
Comparing (4.21) and (4.24) we get
\begin{equation} % 4.25
\hat{u}(t)=\frac{1}{\beta_t} E_Q \Big[ \Big( D_t \hat{X}(T)-\hat{X}(T)\int_t^T D_t
\Big( \frac{\alpha_s}{\beta_s}\Big) d\tilde{B}(s)\Big) \mid \cf_t\Big].
\end{equation}
Using Bayes' rule we conclude

\begin{theo}
Suppose $\hat{u}\in\ca_{\cf}$ is an optimal portfolio for the problem
\begin{equation} % 4.26
\sup_{u\in \ca_{\cf}} E[U(X^{(u)}(T),\omega)]
\end{equation}
with
$$
dX^{(u)}(t)=\alpha_t u(t)dt+\beta_t u(t)dB(t).
$$
Then
\begin{equation} % 4.27
\hat{u}(t)=\frac{1}{\beta_t N_t} E \Big[ N_T \Big( D_t\hat{X}(T)-\hat{X}(T) \int_t^T D_t
\Big( \frac{\alpha_s}{\beta_s}\Big) d\tilde{B}(s)\Big) \mid \cf_t\Big]
\end{equation}
and
\begin{equation}
\hat{X}(T)=(U')^{-1}(M_T)\,,
\end{equation}
where $M_t$ is given by (4.19).
\end{theo}

\begin{coro}
Suppose
$$
U(x)=\frac{1}{\gamma} x^\gamma F(\omega)
$$
for some $\cf_T$-measurable bounded $F$.
Then
\begin{equation} % 4.29
\hat{u}(t)=\hat{X}(t)\frac{1}{\beta_t}
\frac{E[N_T(D_tY-Y\int_t^T D_t(\frac{\alpha_s}{\beta_s})d\tilde{B}(s))\mid \cf_t]}
{E[N_TY \mid \cf_t]}
\end{equation}
where
\begin{equation} % 4.30
Y=\Big[ \frac{1}{F}\exp \Big( -\int_0^T \frac{\alpha_s}{\beta_s}dB(s)
-\tfrac{1}{2}\int_0^T \Big(\frac{\alpha_s}{\beta_s}\Big)^2 ds\Big)\Big]^{\frac{1}{\gamma-1}}.
\end{equation}
\end{coro}

\noindent
{\sc Proof.} \
In this case we get
$$
\hat{X}(T)=\Big( \frac{M_T}{F}\Big)^{\frac{1}{\gamma-1}}=M_0^{\frac{1}{\gamma-1}}Y
$$
and
$$
\hat{X}(t)=E_Q[\hat{X}(T) \mid \cf_t]=
\frac{M_0^{\frac{1}{\gamma-1}} E[N_TY\mid \cf_t]}{N_t}\,.
$$
Therefore the result follows from (4.27).

\hfill

\eof

%\section{REFERENCES}

\end{document}